\documentstyle[12pt,amssymb]{amsart}

\newtheorem{theorem}{Theorem}[section]
\newtheorem{defn}[theorem]{Definition}

\newtheorem{lemma}[theorem]{Lemma}

\newtheorem{fact}[theorem]{Fact}

\newtheorem{eple}[theorem]{Example}
\newtheorem{rmk}[theorem]{Remarks}
\newtheorem{dsc}[theorem]{Discussion}
\newtheorem{nota}[theorem]{Notation}

\newsavebox{\indbin}
\savebox{\indbin}{\begin{picture}(0,0)
\newlength{\gnu}
\settowidth{\gnu}{$\smile$} \setlength{\unitlength}{.5\gnu}
\put(-1,-.65){$\smile$} \put(-.25,.1){$|$}
\end{picture}}
\newcommand{\ind}[3]
{\mbox{$\begin{array}{ccc}
\mbox{$#1$} & \usebox{\indbin} & \mbox{$#2$} \\
        & \mbox{$#3$} &
\end{array}$}}

\newcommand{\nind}[3]
{\mbox{$\begin{array}{ccc}
\mbox{$#1$} & \mbox{$\not \, \downarrow$} & \mbox{$#2$} \\
        & \mbox{$#3$} &
\end{array}$}}
\newcommand{\be}{\begin{enumerate}}
\newcommand{\bd}{\begin{defn}}
\newcommand{\bt}{\begin{theorem}}
\newcommand{\bl}{\begin{lemma}}
\newcommand{\ee}{\end{enumerate}}
\newcommand{\ed}{\end{defn}}
\newcommand{\et}{\end{theorem}}
\newcommand{\el}{\end{lemma}}

\begin{document}
\title{The geometry of linear regular types}
\author{Tristram de Piro}
\address{Mathematics Department \\ University of Edinburgh\\  \\
Edinburgh, EH9 3JZ} \email{depiro@@maths.ed.ac.uk}
\thanks{The author was supported by the William Gordon Seggie Brown Research Fellowship}
\begin{abstract}
This paper is concerned with extending the results of \cite{dePK}
in the context of the solution set $D$ of a regular $Lstp$ defined
over $\emptyset$ in a simple theory $T$. In \cite{P}, a notion of
$p$-weight is developed for regular types in stable theories. Here
we show that the corresponding notion holds in simple theories and
give a geometric analysis of associated structures $G(D)$ and
$G(D)^{large}$, the former of which appears in \cite{dePK}. We
show that $D$ is linear iff $G(D)$ and $G(D)^{large}$ (localized,
resp) are both modular with respect to the $p$-closure operator
$cl_{p}$. Finally, we show that modularity of $G(D)^{large}$
provides a local analogue of $1$-basedness for the theory $T$.
\end{abstract}

\maketitle

\section{Preliminaries}

For convenience we will assume that the ambient theory $T$ is supersimple. In \cite{K1}, Kim shows that forking inside simple theories satisfies the Independence Theorem over a model ${\mathcal M}$.\\

In order to apply the independence theorem over parameters, the notion of Lascar strong type is introduced. As is shown in $\cite{BPW}$, if $T$ is supersimple, then $T$ has elimination of hyperimaginaries and the notion of Lascar Strong Type simplifies to the following;\\

$Lstp(\bar a/A)=Lstp(\bar b/A)$ iff $tp(\bar a/acl(A))=tp(\bar b/acl(A))$ iff $stp(\bar a/A)=stp(\bar b/A)$\\

where $acl(A)$ denotes the algebraic closure of $A$ \emph{in ${\mathcal M}^{eq}$}.\\

We can then apply the Independence Theorem for Lascar strong types;\\

If $\ind{\bar c}{\bar d}{A}, Lstp(\bar a/A)=Lstp(\bar b/A)$ and $\ind{\bar a}{\bar c}{A}, \ind{\bar b}{\bar d}{A}$\\

then the nonforking extensions of $tp(\bar a/A)=tp(\bar b/A)$ to $A\bar c$ and $A\bar d$ respectively can be amalgamated.\\

Given a complete type $p(\bar x,b)$, where $b$ denotes a possibly infinite set of parameters, we define the parallelism class of $p(\bar x,b)$ to be;\\

${\mathfrak B}=\{p(\bar x,c):E(c,b)\}$\\

where $E$ is the transitive closure of the relation\\

$R(c,b)\equiv p(\bar x,c)\cup p(\bar x,b)$ have a common non forking extension\\

As is shown in $\cite{KP}$, if $b$ is an amalgamation base, then $E$ is a type definable equivalence relation on $tp(b)$. If $T$ is supersimple, then the parameter set of $stp(a/B)$ is an amalgamation base and we define the canonical base $C=Cb(Lstp(a/B))$ to be the $E$ class of $acl(B)$. Again, assuming $T$ supersimple, $E$ is the intersection of definable equivalence relations on $tp(acl(B))$ and we may take $C$ to be a possibly infinite set of parameters in ${\mathcal M}^{eq}$. Note the assumption that $T$ is supersimple is not critical in what follows provided we work with ${\mathcal M}^{heq}$ instead of ${\mathcal M}^{eq}$. In general, we do not assume that our parameter sets are algebraically closed. As the notion of non-forking is invariant under algebraic closure in ${\mathcal M}^{eq}$, we often implicitly replace a parameter set $B$ by its algebraic closure in ${\mathcal M}^{eq}$, hoping this will not cause confusion.\\

We will require the following facts about canonical bases as given in \cite{KP}, \cite{K2} and \cite{Wag};\\

\begin{fact}

1. The Independence Theorem holds for the restriction of a Lstp
over $A$ to
the base $C\subset acl(A)$. \\

Let $A\subset B$ be sets and $\bar a$ a tuple, then;\\

2. $\ind{\bar a}{B}{A}$ iff $Cb(Lstp(\bar a/B))\subset acl(A)$. As a consequence, if $C=Cb(Lstp(\bar a/A))$, then $\ind{\bar a}{A}{C}$ and $\ind{\bar a}{C}{A}$\\

3. If $D=Cb(Lstp(\bar a/B))$ and $\ind{\bar a}{B}{A}$, then, using the fact that $C$ and $D$ are amalgamation bases, $dcl(C)=dcl(D)$.\\

4. If $\{\bar a_{i}:i<\omega\}$ is a Morley sequence in $Lstp(\bar a/A)$, then $C=Cb(Lstp(\bar a/A))\subset dcl(\bar a_{i}:i<\omega)$\\

\end{fact}

\begin{defn}
A pregeometry is a set $S$ with a closure operation $cl:P(S)\rightarrow P(S)$ satisfying the following axioms found in \cite{P};\\

1. If $A\subseteq S,$ then $A\subseteq cl(A), cl(A)=cl(cl(A))$.\\

2. If $A\subseteq B\subseteq S$, then $cl(A)\subseteq cl(B)$.\\

3. If $A\subseteq S$, $a,b\in S$, then $a\in cl(Ab)\setminus cl(A)$ implies $b\in cl(Aa)$.\\

4. If $a\in S$ and $a\in cl(A)$, then there is some finite $A_{0}\subset A$ with $a\in cl(A_{0})$.\\

We say that $(S,cl)$ is modular if for $A,B$ finite dimensional
closed subsets of $S$, $dim(A\cup B)=dim(A)+dim(B)-dim(A\cap B)$.

\end{defn}

\begin{rmk}

A necessary and sufficient condition for modularity of a pregeometry $(S,cl)$ is the following;\\

Whenever $a,b\in S$, $B\subset S$ is closed and finite
dimensional, $dim(ab)=2$ and $dim(ab/B)\leq 1$, then there is
$c\in cl(ab)\cap cl(B)$ with $c\notin cl(\emptyset)$ (*)

\end{rmk}

\section{Regular Types and p-weight}

Let $p$ be a non-algebraic complete Lascar strong type over
$\emptyset$. Recalling the definition of orthogonality in simple
theories, see for example \cite{Wag}, we say that $p$ is
\emph{regular} if it is orthogonal to all its forking extensions.

\begin{lemma}

If $p$ is regular, the realisations $D$ of $p$ form a pregeometry with the the closure operation $cl$ given by $cl(A)=\{x\in p:\nind{x}{A}{\emptyset}\}$.\\

\end{lemma}

\begin{pf}

We check the axioms, 2 is trivial and 4 follows from the finite
character of forking. 3 follows immediately from forking symmetry
and all the work is in showing that 1 holds, namely we have to see
that if $A\subset p$, $a,b_{1}\ldots b_{n}$ is a tuple in $p$ such
that $\nind{b_{i}}{A}{\emptyset}$ for each $i$ and
$\nind{a}{b_{1},\ldots,b_{n}}{\emptyset}$ then in fact
$\nind{a}{A}{\emptyset}$. Suppose not, so $a$ realises a non
forking extension of $\rho$ to $A$. Each $b_{i}$ realises a
forking extension of $\rho$ to $A$ so by definition of regularity,
we must have that $\ind{a}{b_{1}}{A}$. Now we just repeat the
argument with $Ab_{1}$ replacing $A$, clearly
$\nind{b_{i}}{Ab_{1}}{\emptyset}$ for $i\geq 2$ and again using
regularity $\ind{a}{b_{2}}{Ab_{1}}$, so we get
$\ind{a}{b_{1}b_{2}}{A}$. After $n$ steps, using transitivity, we
have that $\ind{a}{b_{1}\ldots b_{n}}{A}$ and so, as
$\ind{a}{A}{\emptyset}$ we get $\ind{a}{b_{1}\ldots
b_{n}}{\emptyset}$. This contradicts the
original hypothesis.\\

\end{pf}

Given a set of parameters $A$, we let $D_{A}=\{x\in
D:\ind{x}{A}{\emptyset}\}$ and define a closure operation $cl_{A}$
by $cl_{A}(B)=\{x\in D_{A}:\nind{x}{B}{A}\}$ for $B\subset D_{A}$.
Imitating the above proof, one easily checks that $(D_{A},cl_{A})$
forms a pregeometry which we refer to as the localisation $p_{A}$
of $p$ to $A$.

Given any pregeometry $(S, cl)$, we use the standard notation $(S',cl')$, as in \cite{dePK}, for the associated geometry. For closed $B\subset S$, we have a notion of dimension $dim(B)$. For closed sets $B\subset C\subset S$, we define $dim(C/B)=dim(C)-dim(B)$ and for arbitrary sets $B,C\subset S$, we define $dim(C/B)=dim(cl(C\cup B)/cl(B))$. Recall that this notion is additive and the same holds for the localised analogue $dim_{A}$.\\

Let $p_{1}$ and $p_{2}$ be types over possibly different sets. We say that $p_{1}$ is \emph{hereditarily orthogonal} to $p_{2}$ if every extension of $p_{1}$ is orthogonal to $p_{2}$. Now fix a regular complete Lstp $p$ defined over $\emptyset$ and define a Lstp $q$ over a domain $B$ to be $p$-simple if, for all $a\in q$, there exists $F$ with $B\subset F$, $\ind{F}{a}{B}$ and $I\subset p_{F}$  such that $tp(a/FI)$ is hereditarily orthogonal to $p$. Given $F$ as above, we say that $F$ \emph{witnesses the} $p$-\emph{simplicity of} $q$. Note that by an automorphism argument and the fact that $p$ was assumed to be defined over $\emptyset$, in order to check $p$-simplicity of a given type $q$, it is sufficient to verify it for any element realising $q$\\

We define a notion of $p$-weight for Lascar strong types $q$ as follows;\\

$w_{p}(q)=min\{\kappa:$ there is $F\supset B$, $I\subset p_{F}$ as above with $dim_{F}(I)=\kappa \}$. \\

$w_{p}(a/B)=w_{p}(Lstp(a/B))$\\

We say that $F$ \emph{witnesses the} $p$-\emph{weight of} $q$ if $F$ witnesses the $p$-simplicity of $q$ and there exists $I$ as above with $dim_{F}(I)=w_{p}(q)$.\\

\begin{lemma}
If $F$ witnesses the $p$-weight of $Lstp(a/B)$ and $G\supset F$
with $\ind{G}{a}{F}$ then $G$ also witnesses the $p$-weight of
$Lstp(a/B)$.
\end{lemma}
\begin{pf}
We clearly have that $\ind{G}{a}{B}$. Let $I\subset p_{F}$ be
independent with $Card(I)=w_{p}(Lstp(a/B))$ and $Lstp(a/FI)$
hereditarily orthogonal to $p$. Then still $Lstp(a/GI)$ is
hereditarily orthogonal to $p$. Moreover,$I\subset p_{G}$ as if
$c\in I$ with $\nind{c}{G}{F}$, then, by the finite character of
forking, we can find a tuple $\bar g$ with $\ind{\bar g}{a}{F}$
and $\nind{c}{\bar g}{F}$, $\nind{c}{a}{F}$. As $Lstp(c/F)$ is
regular, this contradicts the fact that $Lstp(c/F)$ has ordinary
weight $1$; a proof of this last fact can be found in for example
\cite{Wag}. Now, by definition of $p$-weight, we must have that
$dim_{G}(I)=w_{p}(Lstp(a/B))$.

\end{pf}

\begin{lemma}

Suppose that q=$Lstp(a/B)$ is $p$-simple and $F$ witnesses the
$p$-weight of $q$. Let $Y=\{b:b\in p_{F},\nind{b}{a}{F}\}$. Then
$Lstp(a/FY)$ is hereditarily orthogonal to $p$ and
$w_{p}(q)=dim_{F}(Y)$.

\end{lemma}

\begin{pf}

By the choice of $F$ and regularity of $p$, we can find
independent $I\subset p_{F}$ such that $Lstp(a/FI)$ is
hereditarily orthogonal to $p$ and $w_{p}(q)=Card(I)$. We can
assume that $\nind{b}{a}{F}$ for all $b\in I$. Otherwise, choose
$I_{0}\subset I$ maximal with $\ind{a}{I_{0}}{F}$. Now replace $F$
by $FI_{0}$ and $I$ by $I\setminus I_{0}$, then the pair
$(FI_{0},I-I_{0})$ work in the definition of $w_{p}(q)$, which is
a contradiction unless $I_{0}=\emptyset$. Hence $I\subset Y$.
Moreover, $I$ is a basis for $Y$ over $F$. If not, we can find
$y\in Y$ with $\ind{y}{I}{F}$ and $\nind{y}{a}{F}$. If
$\ind{y}{a}{FI}$, then $\ind{y}{aI}{F}$, contradicting the
definition of $Y$. Therefore, $\nind{y}{a}{FI}$, contradicting the
fact that $Lstp(a/FI)$ is hereditarily orthogonal to $p$.
\end{pf}

\begin{theorem}

We have the following properties of $p$-weight, $0$ only holding with the assumption that $T$ is supersimple;\\

0. $w_{p}$ is finite: If $Lstp(a/B)$ is $p$-simple, then $w_{p}(a/B)$ is finite.\\

1. Non-Forking: If $B\subset C$, $\ind{a}{C}{B}$ then $Lstp(a/B)$ is $p$-simple iff $Lstp(a/C)$ is $p$-simple and $w_{p}(a/B)=w_{p}(a/C)$.\\

2. Extension: If $Lstp(a/B)$ is $p$-simple and $B\subset C$ then $Lstp(a/C)$ is $p$-simple and $w_{p}(a/C)\leq w_{p}(a/B)$\\

3. Additivity: If $Lstp(a/B)$ and $Lstp(b/B)$ are $p$-simple, then so is $Lstp(ab/B)$ and $w_{p}(ab/B)=w_{p}(a/Bb)+w_{p}(b/B)$. \\

4. Algebraicity: If $Lstp(a/B)$ is $p$ simple and $b\in acl(aB)$ then $Lstp(b/B)$ is $p$-simple and $w_{p}(b/B)\leq w_{p}(a/B)$\\

5. Finite Character: If $Lstp(a/B)$ is $p$-simple and $B\subset C$ then there exists a finite $\bar c\subset C$ such that $w_{p}(a/C)=w_{p}(a/B\bar c)$\\

6. Permutation: $Lstp(ab/B)$ is $p$-simple iff $Lstp(ba/B)$ is
$p$-simple and $w_{p}(ab/B)=w_{p}(ba/B)$
\end{theorem}

\begin{pf}

These are adaptations of the corresponding properties in the stable case;\\

0. Let $F\supset B$ witness $p$-weight for $Lstp(a/B)$ and $Y$ the set given by Lemma 2.3. Then, it's sufficient to prove that $dim_{F}(Y)$ is finite. Suppose not, then we can find an infinite sequence $\{c_{i}:0\leq i<\omega\}$ independent over $F$ such that $\nind{c_{i}}{a}{F}$. By transitivity of non-forking, we have that $\nind{c_{i}}{a}{Fc_{0}\ldots c_{i-1}}$ for each $i$. This gives an infinite forking chain and contradicts the fact that $T$ is supersimple.\\

1. Suppose $Lstp(a/B)$ is $p$-simple and let $F$ witness the
$p$-simplicity. Choose $F'\equiv_{aB}F$ with $\ind{F'}{C}{aB}$.
Then $\ind{F'}{aC}{B}$ and $\ind{F'}{a}{C}$. By automorphism, $F'$
witnesses the $p$-simplicity of $Lstp(a/B)$, so we can find
$I'\subset p_{F'}$ with $tp(a/F'I')$ hereditarily orthogonal to
$p$. Then $tp(a/F'I'C)$ is hereditarily orthogonal to $p$ which
shows that $Lstp(a/C)$ is $p$-simple. Now suppose that $F$
witnesses the  $p$-simplicity of $Lstp(a/B)$ and $I\subset p_{F}$
with $dim_{F}(I)=\kappa$, $Lstp(a/FI)$ hereditarily orthogonal to
$p$. By the same argument, we may assume that $\ind{F}{a}{C}$.
Then $tp(a/FCI)$ is hereditarily orthogonal to $p$ and
$dim_{FC}(I)\leq\kappa$, which shows that $w_{p}(a/C)\leq
w_{p}(a/B)$. Conversely, suppose that $F$ witnesses $p$-simplicity
of $Lstp(a/C)$ and $I\subset p_{F}$ with $dim_{F}(I)=\kappa$,
$Lstp(a/FI)$ hereditarily orthogonal to $p$. Then $\ind{a}{F}{B}$
as $Lstp(a/C)$ is a non-forking extension of $L
stp(a/B)$. This shows that $Lstp(a/B)$ is $p$-simple and $w_{p}(a/B)\leq w_{p}(a/C)$.\\

2. Let $B_{1}$ witness $p$-weight for $Lstp(a/B)$. By the usual
arguments, we may assume that $\ind{B_{1}}{aC}{B}$, hence
$\ind{a}{B_{1}}{C}$. Now choose $I_{1}\subset p_{B_{1}}$ with
$w_{p}(Lstp(a/B))=dim_{B_{1}}(I_{1})$ and $Lstp(a/B_{1}I_{1})$
hereditarily orthogonal to $p$. Then still $Lstp(a/B_{1}CI_{1})$
$(*)$ is hereditarily orthogonal to $p$. Now choose $I_{1}'\subset
I_{1}$ with $I_{1}'\subset p_{B_{1}C}$ and let $J_{1}$ be a basis
for $I_{1}'$ over $B_{1}C$. Then we claim that
$Lstp(a/B_{1}CJ_{1})$ is hereditarily orthogonal to $p$. As
$dim_{B_{1}C}(J_{1})=dim_{B_{1}C}(I_{1}')\leq
dim_{B_{1}}(I_{1})=w_{p}(Lstp(a/B))$, this is sufficient to prove
2. In order to show the claim, let $F\supset B_{1}CJ_{1}$ and
suppose that $d\in p_{F}$. Without loss of generality, we may
assume that still $I_{1}'\subset p_{F}$ and $J_{1}$ is still a
basis for $I_{1}'$ over $F$, otherwise take the corresponding
subsets. By regularity of $p$, $\ind{b}{I_{1}\setminus
I_{1}'}{F}$, hence $b
\in p_{F, I_{1}\setminus I_{1}'}$. If $\nind{b}{I_{1}'}{F, I_{1}\setminus I_{1}'}$, then $b\in cl_{F, I_{1}\setminus I_{1}'}(I_{1}')$ and $I_{1}'\subset cl_{F, I_{1}\setminus I_{1}'}(J_{1})$ implies $b\in cl_{F, I_{1}\setminus I_{1}'}(J_{1})$ by transitivity of $cl_{F,I_{1}\setminus I_{1}'}$. Therefore, $\nind{b}{J_{1}}{F, I_{1}\setminus I_{1}'}$ contradicting the fact that $J_{1}\subset F$ and we conclude that $b\in p_{FI_{1}}$. By $(*)$, we have that $\ind{b}{a}{FI_{1}}$ and the previous argument shows $\ind{b}{I_{1}}{F}$, hence $\ind{b}{a}{F}$ as required.\\

3. For the first part of $3$, suppose that $Lstp(a/B)$ and
$Lstp(b/B)$ are $p$-simple and let $B_{1},B_{2}\supset B$ witness
their $p$-simplicity. Choose $b'\equiv_{Ba}b$, with
$\ind{b'}{B_{1}}{Ba}$ and find $C\equiv_{B}B_{2}$ with $C$
witnessing $p$-simplicity for $Lstp(b'/B)$. Choose
$B_{2}'\equiv_{Bb'}C$ with $\ind{B_{2}'}{aB_{1}}{Bb'}$, then, as
$\ind{B_{2}'}{b'}{B}$, we have that $\ind{B_{2}'}{ab'B_{1}}{B}$.
Therefore, $\ind{B_{2}'}{b'}{B_{1}a}$ and
$\ind{B_{2}'}{a}{B_{1}}$. This gives $ab'\equiv_{B}ab$ and
$\ind{ab'}{B_{1}B_{2}'}{B}$ with $B_{1}$ witnessing $p$-simplicity
for $Lstp(a/B)$ and $B_{2}'$ witnessing $p$-simplicity for
$Lstp(b'/B)$. Then we can find $I_{1},I_{2}$ such that
$Lstp(a/B_{1}I_{1})$ and $Lstp(b'/B_{2}'I_{2})$ are both
hereditarily orthogonal to $p$. We then claim that
$Lstp(ab'/B_{1}B_{2}'I_{1}I_{2})$ is hereditarily orthogonal to
$p$. If not, then we can find $F\supset B_{1}B_{2}'I_{1}I_{2}$ and
$c\in p_{F}$ such that $\nind{ab'}{c}{F}$. Then $
\nind{c}{b'}{F}$ or $\nind{c}{a}{Fb'}$. In either case, we have a contradiction. Hence, $Lstp(ab'/B)=Lstp(ab/B)$ is $p$-simple. \\

In order to show that $w_{p}$ is additive, let $B_{1}$ witness
$p$-weight for $Lstp(a/B)$ with corresponding $I_{1}$. By the
first part, $Lsp(ab/B)$ is $p$-simple so we can find $B_{2}$
witnessing the $p$-weight of $Lstp(ab/B)$. Similar to the above,
we may assume that $\ind{B_{2}}{abI_{1}B_{1}}{B}$ and
$\ind{ab}{B_{1}B_{2}}{B}$. By Lemma 2.2, we may assume that
$B_{1}B_{2}$ witness $p$-weight both for $Lstp(a/B)$ and
$Lstp(ab/B)$. Now choose $I_{2}$ maximal with $I_{1}I_{2}$
independent over $B_{1}B_{2}$ and such that
$\nind{c}{ab}{B_{1}B_{2}}$ for $c\in I_{2}$. Using Lemma 2.3, we
have that $Lstp(ab/B_{1}B_{2}I_{1}I_{2})$ is hereditarily
orthogonal to $p$ and
$w_{p}(ab/B)=dim_{B_{1}B_{2}}(I_{1}I_{2})=dim_{B_{1}B_{2}}(I_{1})+dim_{B_{1}B_{2}}(I_{2})=Card(I_{2})+w_{p}(a/B)$,
the last equality following form the fact that $B_{1}B_{2}$
witnesses $p$-weight for $Lstp(a/B)$. It is therefore sufficient
to prove that $w_{p}(b/aB)=Card(I_{2})$. Choose parameters
$B_{3}\supset aB$ witnessing $p$-weight for $Lstp(b/aB)$. We may
assume that $\ind{b}{B_{1}B_{2}B_{3}I_{1}}{aB}$ and
$\ind{B_{3}}{bB_{1}B_{2}I_{1}I_{2}}{aB}$. By Lemma 2.2,
$B_{1}B_{2}B_{3}I_{1}$ witnesses the $p$-weight of $Lstp(b/aB)$
$(1)$. We have that $I_{2}$ is independent over $aB_{1}B_{2}I_{1}$
$(*)$, as $\ind{I_{2}}{a}{B_{1}B_{2}I_{1}}$ by the fact that
$Lstp(a/B_{1}I_{1})$ is hereditarily orthogonal to $p$ and $I_{2}$
is independent over $I_{1}B_{1}B_{2}$ by definition. Hence, by
choice of $B_{3}$, $I_{2}$ is independent over
$B_{1}B_{2}B_{3}I_{1}$ $(2)$. Now, if $c\in I_{2}$, then
$\nind{c}{b}{B_{1}B_{2}B_{3}I_{1}}$ $(3)$. Otherwise, by $(*)$,
$\ind{c}{ab}{B_{1}B_{2}I_{1}}$. Now $\ind{c}{I_{1}}{B_{1}B_{2}}$
which gives $\ind{c}{ab}{B_{1}B_{2}}$. This contradicts the choice
of $I_{2}$. Combining $(1),(2),(3)$ and Lemma 1.3, it's sufficient
to show that if $d\in p_{B_{1}B_{2}B_{3}I_{1}I_{2}}$ $(**)$, then
$\ind{d}{b}{B_{1}B_{2}B_{3}I_{1}}$. We already have that
$Lstp(ab/B_{1}B_{2}I_{1! }I_{2})$ is hereditarily orthogonal to
$p$, hence so is $Lstp(ab/B_{1}B_{2}B_{3}I_{1}I_{2})$ and
$Lstp(b/B_{1}B_{2}B_{3}I_{1}I_{2})$.
Now, if $(**)$ holds, we have that $\ind{d}{b}{B_{1}B_{2}B_{3}I_{1}I_{2}}$. Then $\ind{d}{I_{2}}{B_{1}B_{2}B_{3}I_{1}}$ implies $\ind{d}{b}{B_{1}B_{2}B_{3}I_{1}}$ as required.\\

\end{pf}

\begin{rmk}

Note that a trivial consequence of 3. is that for a tuple of
elements $\bar a\subset p$ and parameters $F$, $w_{p}(\bar
a/F)=dim_{F}(\bar a')$ where $\bar a'\subset \bar a$ is the
subtuple of elements in $p_{F}$.

\end{rmk}

4. Suppose $C\supset B$ witnesses $p$-weight for $a$, then, as $b\in acl(aB)$, $\ind{C}{b}{B}$. Choose $I\subset p_{C}$ with $Lstp(a/CI)$ $(*)$ hereditarily orthogonal to $p$ and $w_{p}(a/B)=Card(I)$. Then still $Lstp(b/CI)$ is hereditarily orthogonal to $p$. For suppose not, then there exists $D\supset CI$ and $d\in p_{D}$ with $\nind{d}{b}{D}$, hence $\nind{d}{a}{D}$, contradicting $(*)$. This clearly shows $4$. \\

As a corollary of 4, we have the following;\\

\begin{lemma}

Suppose that $Lstp(a/B)$ is $p$-simple, $B\subset C$ and
$E=Cb(a/C)$, then $Lstp(E/B)$ is $p$-simple with finite
$p$-weight.

\end{lemma}
\begin{pf}
By $2$ of Fact 1.1 and the assumption that $T$ is supersimple, we can find a finite set $E_{0}\subset E$ with $E\subset acl(E_{0})$. By $4$ of Fact 1.1, we can find a finite Morley sequence with $E_{0}\subset dcl(a_{i}:1\leq i\leq n)$, hence $E\subset acl(a_{i}:1\leq i\leq n)$. Now, using 3,4 of Theorem 2.4, we have that $Lstp(E/B)$ is $p$-simple and $w_{p}(E/B)\leq nw_{p}(a/B)$ as required.\\
\end{pf}

We also require the following, the proof is the same as for the
stable case, given in \cite{P};

\begin{lemma}

Suppose $Lstp(a/B)$ is $p$-simple and $B\subset C$ with
$w_{p}(a/B)=w_{p}(a/C)=n$. Then, if $E=Cb(a/C)$, $E\subset
cl_{p}(B)$.

\end{lemma}

\begin{rmk}

We can consider $C$ as a finite tuple in ${\mathcal M}^{eq}$ up to
interalgebraicity.

\end{rmk}

5. Choose $F\supset B$ witnessing $p$-weight for $Lstp(a/B)$ with
$\ind{F}{aC}{B}$. Then we can find $I\subset p_{F}$ with
$Card(I)=w_{p}(a/B)$ and $w_{p}(a/FI)=0$. Using 3, we also
calculate $w_{p}(I/aF)=0$. Now, by 1, $w_{p}(a/C)=w_{p}(a/FC)$.
Using 2,3, $w_{p}(a/FC)=w_{p}(aI/FC)=w_{p}(I/FC)$. Hence,
$w_{p}(a/C)=w_{p}(I/FC)$. Now, by the remarks after 3, we can find
$J\subset I'\subset I$ with $I'$ the subtuple belonging to
$p_{FC}$, $J$ a basis for $I'$ over $FC$ and
$w_{p}(I/FC)=w_{p}(J/FC)$. Now, if $i\in I\setminus I'$, then
$\nind{i}{FC}{\emptyset}$ and if $i\in I'\setminus J$, then
$\nind{i}{J}{FC}$, hence, if $i\in I\setminus J$, then
$\nind{i}{FCJ}{\emptyset}$ $(*)$. By the finite character of
forking, we can find a finite tuple $\bar c\subset C$ such that
$(*)$ holds replacing $FCJ$ by $F\bar c J$. Then still $I'$ is the
subtuple of $I$ belonging to $p_{F\bar c}$ and $J$ is a basis for
$I'$ over $F\bar c$. Hence $w_{p}(I/FC)=w_{p}(I/F\bar c)$. Using
2,3 again, $w_
{p}(a/F\bar c)=w_{p}(I/F\bar c)$ and by 1 again, $w_{p}(a/F\bar c)=w_{p}(a/B\bar c)$ as required.\\

6. The proof of $6$ is trivial by the definitions of $p$-simplicity, $p$-weight and the fact that for $B\subset C$, $\ind{ab}{C}{B}$ iff $\ind{ba}{C}{B}$.\\

In order to make the definitions of Section 3, we require one more notion.\\

\begin{defn}

$Lstp(a/B)$ is $p$-pure if for every $B\subset C$,
$w_{p}(a/B)=w_{p}(a/C)$ iff $\ind{a}{C}{B}$.

\end{defn}

The fundamental results on $p$-pure types are the following;\\

\begin{lemma}

Suppose that $Lstp(a/X)$ is $p$-simple and $w_{p}(e/X)\geq 1$ for
all $e\in acl(aX)\setminus acl(X)$ with $Lstp(e/X)$ $p$-simple.
Then $Lstp(a/X)$ is $p$-pure.
\end{lemma}

\begin{pf}

Choose $F\supset X$ witnessing $p$-weight for $Lstp(a/X)$ and
$I\subset p_{F}$ with $Lstp(a/FI)$ hereditarily orthogonal to $p$
and $dim_{F}(I)=w_{p}(a/F)$. We first show that $a$ and $I$ are
equidominant over $F$, see \cite{Wag} for the relevant definition.
Suppose that $\nind{c}{I}{F}$ $(*)$. As
$w_{p}(a/FI)=w_{p}(I/Fa)=0$, we have $w_{p}(a/Fc)=w_{p}(I/Fc)$ by
2,3 of Theorem 2.4. By $(*)$ and regularity of $p$,
$w_{p}(I/Fc)<w_{p}(I/F)$. Hence, $w_{p}(a/Fc)<w_{p}(a/F)$, which,
by 1 of Theorem 2.4, implies $\nind{c}{a}{F}$. This shows that $a$
dominates $I$ over $F$. We want to see that $I$ dominates $a$ over
$F$ as well. Suppose not, then we find $c$ with $\ind{c}{I}{F}$
and $\nind{c}{a}{F}$. Let $E=Cb(Lstp(aI/cF))$, then, by
$\nind{aI}{c}{F}$ and 2. of Fact 1.1,  $E\subset acl(cF)$ and we
can find $e\in E\setminus acl(F)$ with $\nind{e}{a}{F}$ $(*)$. By
4 of Fact 1.1, we can find a Morley sequence $\{a_{1}I_{1},\ldots,
a_{n}I_{n}\}$ realising $Lstp(aI/cF)$ such that $e\in
acl(a_{1}I_{1},\ldots,a_{n}I_{n})$. Then, by 3,4 of Theorem 2.4,
$Lstp(e/FI_{1}! \ldots I_{n})$ is $p$-simple and
$w_{p}(e/FI_{1}\ldots I_{n})=0$. As $\ind{I_{1}\ldots
I_{n}}{c}{F}$ and $e\in acl(cF)$, we have, by 1 of Theorem 2.4,
that $Lstp(e/F)$ is $p$-simple and $w_{p}(e/F)=0$. Finally, let
$C'=Cb(Lstp(eF/aX))$, then, by $(*)$ and 2. of Fact 1.1, we can
find $c'\in C'$ with $c'\in acl(aX)\setminus acl(X)$. As
$Lstp(a/X)$ is $p$-simple, by 4 of Theorem 2.4 we have that
$Lstp(c'/X)$ is $p$-simple. Again, by 4 of Fact 1.1, we can find a
Morley sequence $\{e_{1}F_{1},\ldots, e_{n}F_{n}\}$ realising
$Lstp(eF/aX)$ such that $c'\in acl(e_{1}F_{1},\ldots,
e_{n}F_{n})$. Then again by 3,4 of Theorem 2.4,
$w_{p}(c'/XF_{1}\ldots F_{n})=0$ and, as $\ind{c'}{F_{1}\ldots
F_{n}}{X}$, by 1 of Theorem 2.4, $Lstp(c'/X)$ is $p$-simple and
$w_{p}(c'/X)=0$ as well. This contradicts the assumption of the
Lemma, so $a$ and $I$ are equidominant over $F$. Finally, suppose
that $X\subset Y$ and $w_{p}(a/X)=w_{p}(a/Y)$. Choose $F$ as above
with $\ind{F}{aY}{X}$ and $I\subset p! _{F}$ such that $a$ and $I$
are equidominant over $F$. Then, by 1 of Theorem 1.4,
$w_{p}(I/FY)=w_{p}(a/FY)=w_{p}(a/Y)=w_{p}(a/X)=w_{p}(I/F)$. Hence,
$I\subset p_{FY}$ and $\ind{I}{Y}{F}$ by regularity of $p$. Then
$\ind{a}{Y}{F}$ as $I$ dominates $a$ over $F$ and $\ind{a}{Y}{X}$.
This shows that $Lstp(a/X)$ is $p$-pure as required.

\end{pf}

\begin{lemma}

Suppose $Lstp(a/B)$ is $p$-simple and $w_{p}(a/B)=n$. \\

Let $B^{reg}_{a}=\{b\in acl(aB):w_{p}(b/B)=0\}$, so $B\subset B^{reg}_{a}$, then $Lstp(a/B^{reg}_{a})$ is $p$-pure and for any $c$ with $Lstp(c/B)$ $p$-simple, $w_{p}(c/B)=w_{p}(c/B^{reg}_{a})$.\\

\end{lemma}

\begin{pf}

We first claim that if $Lstp(e/B^{reg}_{a})$ is $p$-simple and
$e\in acl(aB^{reg}_{a})\setminus acl(B^{reg}_{a})$ then
$w_{p}(e/B^{reg}_{a})\geq 1$ $(*)$. As $e\in acl(aB)$, by 4 of
Theorem 1.4, $Lstp(e/B)$ is $p$-simple. Suppose $(*)$ fails, then,
using $5$, we can find a finite tuple $\bar b\in B^{reg}_{a}$ with
$w_{p}(e/B\bar b)=0$. Now, using 2,3,6 and the definition of
$B^{reg}_{a}$, we can easily calculate $w_{p}(e/B)=w_{p}(\bar b
e/B)=w_{p}(e\bar b/B)=0$. As $e\in acl(aB)$, we conclude that
$e\in B^{reg}_{a}$, contradicting the assumption. Now by $(*)$ and
Lemma 2.10, we conclude that $Lstp(a/B^{reg}_{a})$ is $p$-pure.
For the second part, again we can find a finite tuple $\bar b\in
B^{reg}_{a}$ with $w_{p}(c/B^{reg}_{a})=w_{p}(c/B\bar b)$. Again,
a simple weight calculation gives that $w_{p}(c/B\bar
b)=w_{p}(c\bar b/B)=w_{p}(c/B)$ as required.

\end{pf}

\begin{lemma}

Suppose that $Lstp(a/B)$ is $p$-simple, $\ind{a}{C}{B}$ and
$Lstp(a/C)$ is $p$-pure, then $Lstp(a/B)$ is $p$-pure.

\end{lemma}

\begin{pf}

Choose $F\supset C$ witnessing the $p$-weight of $Lstp(a/C)$, so
$F$ also witnesses the $p$-weight of $Lstp(a/B)$. If $I\subset
p_{F}$ is chosen with $Lstp(a/FI)$ hereditarily orthogonal to $p$
and $dim_{F}(I)=w_{p}(a/F)$, then, as we have seen, $a$ and $I$
are weight equivalent over $F$. As is easily seen, $Lstp(I/F)$ is
$p$-pure, hence, as $Lstp(a/F)$ is $p$-pure, an easy check gives
that $a$ and $I$ are domination equivalent over $F$. Now repeat
the proof at the end of Lemma 2.10.

\end{pf}

\section{Linearity and $1$-Basedness}

\begin{defn}

We say that $D$ is linear if the following holds;\\

If $ab$ is a pair in $D$ with $Lstp(ab/B)$ $p$-pure having $p$-weight $1$, then $w_{p}(C/\emptyset)\leq 1$ where $C=Cb(Lstp(ab/B))$.\\
\end{defn}

We also introduce the following $2$ objects.\\

\begin{defn}

$G(D)=\{c:Lstp(c/\emptyset)$ is $p$-simple of $p$-weight $1\}$\\

and \\

$G(D)^{large}=\{c:Lstp(c/\emptyset)$ is $p$-simple of finite $p$-weight$\}$\\

with corresponding localised structures;\\

$G(D)_{A}=\{c\in G(D):w_{p}(c/A)=1\}$\\

and \\

$G(D)_{A}^{large}=\{c\in G(D)^{large}:w_{p}(c/A)\ is\ finite\}$\\

\end{defn}

Note that by 2 of Theorem 2.4, $G(D)^{large}_{A}=G(D)^{large}$. We define a closure operator $cl_{p}$ on $G(D)^{large}$ by $cl_{p}(B)=\{c\in G(D)^{large}:w_{p}(c/B)=0$. We also have a corresponding operator $cl_{p}$ by restriction to $G(D)$ and localised operators $cl_{p,A}$ on $G(D)_{A}$ and $G(D)_{A}^{large}$. \\

\begin{lemma}

The following properties hold for $cl_{p}$. \\

1. $cl_{p}$ is transitive for $G(D)$(local.resp.) and $G(D)^{large}$(local.resp.)\\

2. $cl_{p}$ is finite for $G(D)$(local.resp) and $G(D)^{large}$(local.resp.).\\

3. $cl_{p}$ satisfies exchange on $G(D)$(local.resp).\\

4. $G(D)$ and $G(D)_{A}$ form pregeometries under $cl_{p}$ and $cl_{p,A}$\\

\end{lemma}
\begin{pf}

1. For suppose that $\bar a\in cl_{p}(\bar b)$ and $\bar b\in cl_{p}(\bar c)$ then $w_{p}(\bar a\bar b\bar c/\emptyset)=w_{p}(\bar c/\emptyset)$ by 3. of Theorem 2.4 and $w_{p}(\bar b\bar a \bar c/\emptyset)=w_{p}(\bar a\bar c/\emptyset)=w_{p}(\bar a/\bar c)+w_{p}(\bar c/\emptyset)$, so $w_{p}(\bar a/\bar c)=0$ and $\bar a\in cl_{p}(\bar c)$. The localised proof is similar.\\

2. For suppose that $B\subset G(D)^{large}$ or $B\subset G(D)$  and $\bar a\in cl_{p}(B)$ then ,by property $4$ of $w_{p}$, there is a finite $\bar b\subset acl(B)$ such that $\bar a\in cl_{p}(\bar b)$. By transitivity of $cl_{p}$ and the fact that algebraic types have $p$-weight $0$, we can assume that $\bar b\in B$.\\

3. For suppose that $aBc\subset G(D)$ and $a\in cl_{p}(Bc)\setminus cl_{p}(B)$. Replacing $B$ by $B^{reg}_{c}$ and using Lemma 2.11 we may assume that $w_{p}(a/B^{reg}_{c})=1$ and $Lstp(c/B^{reg}_{c})$ is $p$-pure. Then, as $w_{p}(a/B^{reg}_{c}c)=0$, by the extension property we must have that $\nind{a}{c}{B^{reg}_{c}}$. Therefore, as $c\in G(D)$ and $Lstp(c/B^{reg}_{c})$ is $p$-pure, $w_{p}(c/B^{reg}_{c}a)=0$, that is $c\in cl_{p}(B^{reg}a)$. Then, by transitivity of $cl_{p}$ on $G(D)^{large}$, $c\in cl_{p}(Ba)$.\\

4. The axioms 1,3 and 4 for a pregeometry in Definition 1.2 are easily verified from the previous parts of the lemma. 2 is a straightforward consequence of the extension property for $w_{p}$. \\

\end{pf}

\begin{defn}

We say that $G(D)$ is linear if $(ab)$ is a pair from $G(D)$ and $B\subset G(D)$ such that $Lstp(ab/B)$ is $p$-pure with $w_{p}(ab/B)=1$, then $w_{p}(C)\leq 1$ where $C=Cb(Lstp(ab/B))$.\\

\end{defn}

We now prove the following;\\

\begin{lemma}

If $D$ is linear then $G(D)$ is linear.

\end{lemma}

\begin{proof}
The proof will be similar to \cite{dePK}.  Let $(ab)$ be a pair from $G(D)$ with $w_{p}(ab)=2$, the case for $w_{p}(ab)=1$ is easier, and suppose that $B\subset G(D)$ with $w_{p}(ab/B)=1$ and $Lstp(ab/B)$ $p$-pure. Let $F$ witness the $p$-weight of $Lstp(ab/B)$ with $\ind{F}{B}{ab}$. Then we can find $cd\subset p_{F}$ such that $w_{p}(ab/Fcd)=0$, $w_{p}(cd/F)=2$ and $w_{p}(ab/F)=2$. By additivity of $p$-weight we must have $w_{p}(cd/abF)=0$ as well. \\

Claim 1: $w_{p}(cd/FB)=1$.\\

As $w_{p}(ab/B)=1$, $F\downarrow_{B}ab$ and $w_{p}$ is invariant under non forking extension we have $w_{p}(ab/FB)=1$. Then\\

$w_{p}(abcd/FB)=w_{p}(ab/cdFB)+w_{p}(cd/FB)=w_{p}(cd/FB)=$\\

  $w_{p}(cd/abFB)+w_{p}(ab/FB)=0+1=1$\\

giving the claim.\\

Now replace $FB$ by $FB^{reg}_{cd}=\{b\in acl(cdFB):w_{p}(b/FB)=0\}$\\

By Lemma 2.11, $Lstp(cd/FB^{reg}_{cd})$ is $p$-pure with $p$-weight $1$ and, by linearity of $D$, $w_{p}(C)\leq 1$ where $C=Cb(Lstp(cd/FB^{reg}_{cd}))$. Then\\

Claim 2: $w_{p}(cd/cl_{p}(cd)\cap acl(FB^{reg}_{cd}))=1$\\

We have that $cd\downarrow_{C}FB$, hence $w_{p}(cd/C)=1$. Then, by additivity of $w_{p}$ and linearity of $D$, we calculate $w_{p}(C/cd)=w_{p}(C)-1=0$, therefore $C\subset cl_{p}(cd)$ As $C\subset cl_{p}(cd)\cap acl(FB^{reg})$ and $cl_{p}(cd)\cap acl(FB^{reg})\subset acl(FB^{reg}_{cd})$ the claim is shown.\\

Claim 3: $w_{p}(ab/W)=1$, where $W=cl_{p,F}(ab)\cap acl(FB^{reg}_{cd})\cup F$\\

We clearly still have that $w_{p}(cd/cl_{p,F}(cd)\cap acl(FB^{reg})\cup F)=1$. Using additivity, $w_{p}(ab/cl_{p,F}(cd))\cap acl(FB^{reg})\cup F)=1$. By transitivity of $p$-closure, we must have that $cl_{p,F}(cd)=cl_{p,F}(ab)$, hence $w_{p}(ab/W)=1$ as required.\\

Now let $C'=Cb(Lstp(ab/B))$. Then\\

Claim 4: $w_{p}(ab/WC')=1$\\

If not, then as $C'\in acl(B)$, $ab\in cl_{p}(FB^{reg}_{cd})$. Again by transitivity of $p$-closure and the definition of $FB^{reg}_{cd}$ we must have $ab\in cl_{p}(FB)$. Then, as $ab\downarrow_{B}F$, $ab\in cl_{p}(B)$, contradicting the fact $w_{p}(ab/B)=1$ and giving the claim.\\

Now $ab\downarrow_{C'}B$ so still $w_{p}(ab/C')=1$ and moreover, by Lemma 2.11, $Lstp(ab/C')$ is still $p$-pure. Then, by definition of $p$-purity, we must have that $ab\downarrow_{C'}W$ and so $C'\in acl(W)$. Then $C'\subset cl_{p,F}(ab)$ and as $C'\downarrow_{ab}F$, we must have $C'\subset cl_{p}(ab)$. Now, by calculating $w_{p}(abC')$, we have that $w_{p}(C')=1$ as required.\\
\end{proof}

\begin{lemma}

If $G(D)$ is linear then $G(D)$ is modular.\\

\end{lemma}

\begin{proof}
As $G(D)$ forms a pregeometry, it is sufficient to check the
criterion $(*)$ in Remarks 1.3. So choose $x_{1}x_{2}$ in $G(D)$
with $w_{p}(x_{1}x_{2})=2$ and $Y$ closed finite dimensional such
that $w_{p}(x_{1}x_{2}/Y)=1$. By finiteness, we can find $\bar
y\subset Y$ such that $w_{p}(x_{1}x_{2}/\bar y)=1$ and
$cl_{p}(\bar y)=Y$. Replace $\bar y$ by $\bar
y^{reg}_{x_{1}x_{2}}$, so  we can assume that
$Lstp(x_{1}x_{2}/\bar y^{reg}_{x_{1}x_{2}})$ is $p$-pure, though
$\bar y^{reg}_{x_{1}x_{2}}$ may no longer be contained in $G(D)$.
Now, using 3. of Fact 1.1 and Lemma 2.12, we can replace the
parameters $\bar y^{reg}_{x_{1}x_{2}}$ by a Morley sequence
$W\subset G(D)$ such that $Lstp(x_{1}x_{2}/W)$ is still $p$-pure
and the canonical base is preserved. By linearity of $G(D)$, we
have that $C=Cb(Lstp(x_{1}x_{2}/W))=Cb(Lstp(x_{1}x_{2}/\bar
y^{reg}_{x_{1}x_{2}})\subset cl_{p}(x_{1}x_{2})\cap cl_{p}(\bar
y)$. As $w_{p}(C)=1$, we can find a tuple $\bar c\in G(D)$ such
that $\bar c$ and ! $C$ are interalgebraic. Then in fact $\bar c$
witnesses the criterion $(*)$ as required.
\end{proof}

The $2$ lemmas combine to give the following theorem.\\

\begin{theorem}

If $D$ is linear then $G(D)$ is modular.

\end{theorem}

Even though $G(D)^{large}$ is not a pregeometry it still makes
sense to talk of the dimension of a closed set.

\begin{defn}

Given $X,Y\subset G(D)^{large}$ closed\\

$dim(X/Y)=max \{w_{p}(\bar a/Y):\bar a\in X\}$ if this set is bounded.\\

and $dim(X/Y)=\infty$ otherwise.\\

\end{defn}

\begin{defn}

$G(D)^{large}$ is modular if the following holds ;\\

For finite dimensional closed $X,Y\subset G(D)^{large}$ $dim(X/Y)=dim(X/X\cap Y)$.\\

\end{defn}

We now prove the following;\\

\begin{theorem}

If $D$ is linear then $G(D)^{large}$ is modular

\end{theorem}

Here the problem is made more difficuly by the fact that
$G(D)^{large}$ is not a pregeometry.

\begin{pf}
We first reduce the problem to a finite one, as in in general $cl_{p}(X)$ will be a very large set! Suppose $G(D)^{large}$ is not modular, then there exists closed sets $X$ and $Y$ such that $dim(X/Y)<dim(X/X\cap Y)$. Taking $\bar x\in X$ so that $w_{p}(\bar x/Y)$ is maximal, by definition we have that $w_{p}(\bar x/X\cap Y)<w_{p}(\bar x/Y)$. By finiteness, I can find $\bar c\subset X\cap Y$ and $\bar y\subset Y$ such that $w_{p}(\bar x/\bar c)<w_{p}(\bar x/\bar y)$ and, moreover, as weight is preserved on both sides, we can take $\bar c$ and $\bar y$ such $cl_{p}(\bar c)=X\cap Y$ and $cl_{p}(\bar y)=Y$. Therefore, it is sufficient to prove that \\

$w_{p}(\bar x/\bar y)=w_{p}(\bar x/\bar c)$ where $cl_{p}(\bar c)=cl_{p}(\bar x)\cap cl_{p}(\bar y)$ (*)\\

We show $(*)$  by induction on $w_{p}(\bar x/\bar y)$ for $\bar x$ and $\bar y$  finite tuples from $G(D)^{large}$.\\

Base Case. $w_{p}(\bar x/\bar y)=1$.\\

Suppose $w_{p}(\bar x)=n$, then I can find $F\downarrow \bar x\bar
y$ and $z_{1}\ldots z_{n}\in p_{F}$ such that $\bar x$ and
$z_{1}\ldots z_{n}$  are weight equivalent over $F$ (*).  As
before, one checks that $w_{p}(z_{1}\ldots z_{n}/F\bar y)=1$.
Without loss of generality, we can assume that $w_{p}(z_{i}/F\bar
y)=1$ for each $i$. Now, adding parameters $e_{1}\ldots
e_{n}\subset cl_{p}(\emptyset)$, we may assume that $Lstp(z_{i}/
e_{1}\ldots e_{n})$ is $p$-pure for all $i$ and all the conditions
are preserved with $F\bar y e_{1}\ldots e_{n}$ replacing $F\bar
y$. We must have that $w_{p}(z_{1}z_{i}/F\bar y e_{1}\ldots e_{n}
)=1$ for all $i$, hence by linearity of $D$, we can find $c_{i}\in
G(D)$ for $i\geq 2$ with $cl_{p}(c_{i})=cl_{p}(z_{1}z_{i})\cap
cl_{p,F}(\bar y e_{1}\ldots e_{n})$. Clearly, $e_{i}\subset
cl_{p}(c_{i})$, so without loss of generality $e_{i}\subset
c_{i}$. Now $w_{p}(c_{i}/z_{1}z_{i})=0$ and
$w_{p}(c_{i}/z_{1})=1$, otherwise $\nind{z_{1}}{c_{i}}{e_{1}}$
  and $z
_{1}\subset cl_{p,F}(\bar y)$. Hence
$\nind{z_{i}}{c_{i}}{z_{1}e_{i}}$. As $w_{p}(z_{i}/e_{i}z_{1})=1$
and $Lstp(z_{i}/e_{i}z_{1})$ is $p$-pure, we have that
$w_{p}(z_{i}/z_{1}c_{i})=0$ so $z_{i}\subset cl_{p}(z_{1}c_{i})$.
We want to show that $w_{p}(c_{2}\ldots c_{n})=n-1$ from which,
taking $\bar c=c_{2}\ldots c_{n}$, we clearly have that
$w_{p}(z_{1}\ldots z_{n}/\bar c)=w_{p}(z_{1}\ldots z_{n}/F\bar y)$
and $cl_{p}(\bar c)=cl_{p}(z_{1}\ldots z_{n})\cap cl_{p,F}(\bar
y)$. Suppose not, say $c_{n}\subset cl_{p}(c_{2}\ldots c_{n-1})$,
then as $z_{n}\subset cl_{p}(z_{1}c_{n})$ and $c_{2}\ldots
c_{n-1}\subset cl_{p}(z_{1}\ldots z_{n-1})$, we have that
$z_{n}\subset cl_{p}(z_{1}\ldots z_{n-1})$ contradicting the fact
that $z_{1}\ldots z_{n}$ are independent realisations of $p_{F}$.
Now, $cl_{p}(\bar z)\cap cl_{p,F}(\bar y)\subset cl_{p,F}(\bar
z)\cap cl_{p,F}(\bar y)\cup F$. Therefore, using 3. of Theorem 2.4
to check that $w_{p}(z_{1}\ldots z_{n}/F\bar y)=w_{p}(z_{1}\ldots
z_{n}/cl_{p,F}(\bar y)\cup F)$, we have that $w_{p}(\bar z/F\bar
y)=w_{p}(\bar z/cl_{p,F}(\bar z)\cap cl_{p,F}(\bar y)\cup F)$. Let
$W=cl_{p,F}(\bar x)\cap cl_{p,F}(\bar y)\cup F$, then, using
$(*)$, $w_{p}(\bar x/W)=w_{p}(\bar x/\bar y)=1$. Finally, we can
assume that
$Lstp(\bar x/\bar y)$ is $p$-pure and one checks that $w_{p}(\bar x/WC)=1$, where $C\subset G(D)^{large}$ is $Cb(Lstp(\bar x/\bar y))$. As in Lemma 1.18, this forces $C\subset cl_{p}(W)$ and then $C\subset cl_{p,F}(\bar x)$, and then $C\subset cl_{p}(\bar x)$. This gives the result, as $w_{p}(\bar x/\bar y)=w_{p}(\bar x/C)=w_{p}(\bar x/cl_{p}(\bar x)\cap cl_{p}(\bar y))$ and clearly $cl_{p}(C)=cl_{p}(\bar x)\cap cl_{p}(\bar y)$, otherwise we could find $z\in cl_{p}(\bar x)\setminus cl_{p}(C)$ such that $w_{p}(\bar x/C)=w_{p}(\bar x/Cz)$, which contradicts $3$ of Theorem 2.4\\

Induction Step.\\

We now inductively assume the result for $\bar x$ and $\bar y$
with $w_{p}(\bar x/\bar y)=m$ and suppose that $w_{p}(\bar x/\bar
y)=m+1$.  Now again we can find $F\downarrow \bar x\bar y$ and
$z_{1}\ldots z_{n}\in p_{F}$ such that $z_{1}\ldots z_{n}$ is
weight equivalent to $\bar x$ over $F$. Then still $w_{p}(\bar
x/F\bar y)=m+1$ and we may assume $z_{i}\notin cl_{p}(F\bar y)$
for some $i$, otherwise $\bar x\in cl_{p}(F\bar y)$ which is not
the case. Using the fact that $w_{p}(z_{1}/F\bar y)=1$ say, then
by a weight calculation we have that $w_{p}(\bar x/z_{1}F\bar
y)=m$. We now temporarily add $F$ to the language, and take
$p$-closure to include $F$. Then, working in $G(D)^{large}_{F}$,
we have that $w_{p}(\bar x/\bar y)=m+1$ and $w_{p}(\bar
x/z_{1}\bar y)=m$. Applying the induction hypothesis to
$G(D)^{large}_{F}$, we can find $c$ in $G(D)^{large}_{F}$ such
that $cl_{p}(c)=cl_{p}(\bar x)\cap cl_{p}(z_{1}\bar y)$. Then
$w_{p}(cz_{1}/\bar y)=1$ as $c\in cl_{p}(z_{1}\bar y)$
and $z_{1}\notin cl_{p}(\bar y)$. Therefore, we can find $d\in G(D)^{large}_{F}$ such that $cl_{p}(d)=cl_{p}(cz_{1})\cap cl_{p}(\bar y)$ and moreover $w_{p}(d)=w_{p}(cz_{1})-1=w_{p}(c)-1=w_{p}(\bar x)-m-1$. As $cl_{p}(cz_{1})\cap cl_{p}(\bar y)=cl_{p}(\bar x)\cap cl_{p}(\bar y)$, this tells us exactly that $w_{p}(\bar x/F\bar y)=w_{p}(\bar x/Fd)$ where $cl_{p,F}(d)=cl_{p,F}(\bar x)\cap cl_{p,F}(\bar y)$. Now letting $C'=Cb(Lstp(\bar x/\bar y)$ and assuming as usual that $Lstp(\bar x/\bar y)$ is $p$-pure, we have that $w_{p}(\bar x/FdC')=m$ otherwise as $C'\in acl(\bar y)$ then $w_{p}(\bar x/F\bar y)<m$ which is not the case. Hence, by $p$-purity, we have that $C'\in cl_{p}(F\bar x)$ and then as $F\downarrow_{\bar x}C'$, $C'\in cl_{p}(\bar x)$. This proves the result.\\
\end{pf}
So we have,

\begin{theorem}

If $D$ is linear then $G(D)$ and $G(D)^{large}$ are both modular.

\end{theorem}

The following result is an easy adaptation of the proof in the
stable case, given in \cite{P}, p269, $(iii)\rightarrow (iv)$.
Here, observe that the corresponding notation to $D(p,A_{0})$ is
$G(D)^{large}$ as we have taken $A_{0}=\emptyset$ and that Lemma
2.1 there corresponds to Lemma 2.7 in this paper.

\begin{theorem}

If $G(D)^{large}$ is modular, then $D$ is linear.

\end{theorem}

Also, the following result is a straightforward adaptation of Theorem 3.10;\\

\begin{theorem}

If $G(D)$ is modular, then $G(D)^{large}$ is modular.

\end{theorem}

Combining these results gives that

\begin{theorem}

$D$ is linear iff $G(D)$ and $G(D)^{large}$ are modular.

\end{theorem}

\begin{rmk}

The localised analogues of Theorem 3.14 replace $G(D)$ and
$G(D)^{large}$ by $G(D)_{A}$ and $G(D)^{large}_{A}$. For
$G(D)^{large}_{A}$ there is nothing to prove and that $G(D)_{A}$
modular implies $G(D)$ linear is straightforward.

\end{rmk}

Modularity of $G(D)^{large}$ can be seen as a local analogue of $1$-basedness for the theory $T$. More precisely, we say that a simple theory $T$ with elimination of hyperimaginaries is $1$-based if, for any sets $A$ and $B$ in a big model ${\mathcal M}$, we have that $\ind{A}{B}{acl(A)\cap acl(B)}$ where $acl$ is taken in the sense of ${\mathcal M}^{eq}$. This is equivalent to the following condition (*) on canonical bases.\\

For any tuple $\bar a$ and parameters $B\subset {\mathcal M}^{eq}$, then $Cb(\bar a/B)\subset acl(\bar a)$ (*) \\

The proof is fairly straightforward; Suppose that $T$ is $1$-based, then $\ind{\bar a}{B}{acl(\bar a)\cap B}$. By 2. of Fact 1.1, we have that (*) holds. Conversely, suppose that (*) holds, then given $\bar a, B$, by 2. of Fact 1.1 again, we must have that $\ind{\bar a}{B}{acl(a)\cap acl(B)}$. Now, $T$ must be $1$-based by the finite character of forking.\\

If $T$ is a simple $1$-based theory and $D$ denotes the solution set of any regular type, then, if $X,Y$ are $p$-closed subsets of $G(D)^{large}$ (therefore algebraically closed), by $1$-basedness we must have that $\ind{X}{Y}{X\cap Y}$. Using 1. of Theorem 2.4, we get that $w_{p}(X/Y)=w_{p}(X/X\cap Y)$, so $G(D)^{large}$ is modular and, in particular, by Theorem 3.14, $D$ is linear. The converse, in general, is false, there are examples of \emph{stable} theories all of whose regular types are linear but which are not $1$-based, see \cite{P}. However, we can show the following, which should be compared with $(*)$ above;\\

\begin{theorem}

If $D$ is linear, then given $a\in G(D)^{large}$ and parameters
$B$ such that $Lstp(a/B)$ is $p$-pure, we have
$Cb(Lstp(a/B))\subset cl_{p}(a)$, where $cl_{p}$ denotes the
$p$-closure operator on $G(D)^{large}$.

\end{theorem}

\begin{pf}
The proof is a rather immediate consequence of the main Theorem 3.14. Let $C=Cb(Lstp(a/B))$, then $C\subset G(D)^{large}$ by Lemma 2.6. Let $E=cl_{p}(a)\cap cl_{p}(C)$. Then by modularity of $G(D)^{large}$, $w_{p}(a/C)=w_{p}(a/cl_{p}(C))=w_{p}(a/E)=w_{p}(a/EC)$. Now, by $p$-purity of $Lstp(a/C)$, we must have that $\ind{a}{E}{C}$, therefore by 3. of Fact 1.1, $C=Cb(Lstp(a/EC))$. By Lemma 2.7, $C\subset cl_{p}(E)$ and hence $C\subset cl_{p}(a)$ as required.\\
\end{pf}

\end{document}